\documentstyle{amltd2004}
\def\ritem#1{\item[{\rm #1}]}
\begin{document}
\annalsline{156}{2002}
\received{May 15, 2000}
\startingpage{197}
\def\bye{\end{document}}
 \font\tenrm=cmr10

%--------------- Author macros ---------------
%for Bbb in amstex
\catcode`\@=11
\font\twelvemsb=msbm10 scaled 1100
\font\tenmsb=msbm10
%\font\ninemsb=msbm7 scaled 1100%msbm9
\font\ninemsb=msbm10 scaled 800
\newfam\msbfam
\textfont\msbfam=\twelvemsb  \scriptfont\msbfam=\ninemsb
  \scriptscriptfont\msbfam=\ninemsb
\def\msb@{\hexnumber@\msbfam}
\def\Bbb{\relax\ifmmode\let\next\Bbb@\else
 \def\next{\errmessage{Use \string\Bbb\space only in math
mode}}\fi\next}
\def\Bbb@#1{{\Bbb@@{#1}}}
\def\Bbb@@#1{\fam\msbfam#1}
\catcode`\@=12

 \catcode`\@=11
\font\twelveeuf=eufm10 scaled 1100
\font\teneuf=eufm10
\font\nineeuf=eufm7 scaled 1100%eufm9
\newfam\euffam
\textfont\euffam=\twelveeuf  \scriptfont\euffam=\teneuf
  \scriptscriptfont\euffam=\nineeuf
\def\euf@{\hexnumber@\euffam}
\def\frak{\relax\ifmmode\let\next\frak@\else
 \def\next{\errmessage{Use \string\frak\space only in math
mode}}\fi\next}
\def\frak@#1{{\frak@@{#1}}}
\def\frak@@#1{\fam\euffam#1}
\catcode`\@=12
%-------------- Author entries --------------------

%-------------- Article Text--------------------

%\intro %(Optional, Introduction)
\input amssym.def
\input amssym.tex

\def\F{{\cal F}}
\def\B{{\cal B}}
\title{The extremal function\\ associated to intrinsic norms}

\shorttitle{The extremal function} % Shortened version for headline title

 \acknowledgements{Research of the author was supported in part by 
NSERC Grant OGP-0046732.}
\author{Pengfei Guan}
\institution{McMaster University, Hamilton, Ontario, Canada
\\
  {\eightpoint {\it E-mail address\/}: guan@math.mcmaster.ca
}}

  \centerline{\bf Abstract}
\vglue6pt
  Through the study of the degenerate complex Monge-Amp\`ere 
equation, we establish the optimal regularity of the extremal function 
associated to intrinsic norms of Chern-Levine-Nirenberg and Bedford-Taylor. 
We prove a conjecture of Chern-Levine-Nirenberg on the extended
intrinsic norms on complex manifolds and verify Bedford-Taylor's representation
formula for these norms in general.

\section{Introduction}

For every compact complex manifold $M$ with boundary, Chern-Levine-Nirenberg 
defined in \cite{CLN} an intrinsic norm on the homology groups $H_k(M,{\bf R})$
($k=1, \ldots, 2n-1$), as the supremum of $C^2$ plurisubharmonic
functions in a certain class. A similar norm $\tilde N$ was also introduced by
Bedford-Taylor \cite{BT4}. These norms are invariants of complex manifolds
and decreasing under holomorphic mappings. In particular, the characterizations of these norms on $H_{2n-1}(M,{\bf R})$ are very important in the study of
holomorphic mappings. Associated to these norms, there is an extremal function which satisfies the following homogeneous complex Monge-Amp\`ere equation:

\begin{eqnarray}\label{e1.12}
\left\{\begin{array}{ll}
(dd^cu)^n = 0 & \mbox{ in } M^0\\
u|_{\Gamma_1 = 1}\\
u|_{\Gamma_0 = 0},
\end{array}\right.
\end{eqnarray}
where $d^c = i(\bar\partial - \partial)$, $M^0$ is the interior of $M$, and $\Gamma_1$ and $\Gamma_0$ are the
corresponding outer and inner boundaries of $M$ respectively.

Based on the dual principle of the calculus of variations, Chern-Levine-Nirenberg conjectured that the intrinsic norm should
be equal to a minimum of a certain other class of $C^2$ plurisubharmonic functions (see (\ref{e1.200})). 
Under the assumption of $C^2$ regularity of the extremal function $u$ in (\ref{e1.12}), Bedford-Taylor verified the conjecture and
observed the following important representation formula in \cite{BT4}, 
\begin{eqnarray}\label{e1.20}
\tilde N(\{\Gamma_1\})=\int_{\Gamma_1} \left(\frac{\partial u}{dr}\right)^n d^cr \wedge (dd^cr)^{n-1},
\end{eqnarray}
where $\tilde N$ (see (\ref{e1.9}) and (\ref{e1.10})) is the extended Chern-Levine-Nirenberg norm by Bedford-Taylor.
Also, $\Gamma_1 $ is the outer boundary of $M$ and $r$ is a defining function of $\Gamma_1$.
The identity (\ref{e1.20}), together with Fefferman's boundary regularity theorem (\cite{F}), plays a crucial role in the study of holomorphic mappings on annular domains in $\bf C^n$ by Bedford-Burns in \cite{BB}. The idea of using the complex Monge-Amp\`ere equation to study holomorphic mappings was discussed in \cite{K} by Kerzman-Kohn-Nirenberg. 

It is clear that the regularity of the extremal function is the key issue. As the equation (\ref{e1.12})
is a degenerate complex Monge-Amp\`ere equation and a part of the
boundary of the manifold $M$ is pseudoconcave, regularity becomes a difficult
problem. The existence of 
Lipschitz solution was shown in \cite{BT4}, and the uniqueness of the solution of the equation is well known (see e.g.,
\cite{BT1},\cite{K} and \cite{CKNS}).  In some special cases, for example on Reinhardt domains (\cite{BT4}) or a
perturbation of them (\cite{B} and \cite{M}),  the extremal function is smooth. Unfortunately, the solution of
(\ref{e1.12}) is not in $C^2$ in general. In \cite{BF}, Bedford-Forn\ae ss constructed counterexamples with solutions
of exact $C^{1,1}$ regularity.  

In this paper, we will prove the optimal $C^{1,1}$ regularity of the equation (\ref{e1.12}). Moreover, we
will establish the formula (\ref{e1.20}) in general without the $C^2$ regularity
assumption on the extremal function.  And we will also verify the 
Chern-Levine-Nirenberg conjecture for $\tilde N$. We will adapt the subsolution method introduced by B. Guan and J. Spruck \cite{GS} and B. Guan \cite{G} to deal with the degenerate complex
Monge-Amp\`ere equation. One of the crucial steps is the construction
of a smooth
subsolution in Proposition~\ref{subsol}.

\bigskip

Before we state our main
results, we first recall the intrinsic norms defined by Chern-Levine-Nirenberg,
and the extensions of Bedford-Taylor. 

\bigskip

Let $\gamma \in H_*(M, {\Bbb R})$ be a homology class in $M$; define
\begin{eqnarray}\label{e1.1}
N\{\gamma\} = \sup_{u \in \F}\ \inf_{T \in \gamma} |T(d^c u \wedge (dd^cu)^{k-1})|, \quad \mbox{ if } \dim\gamma = 2k-1;
\end{eqnarray}
\begin{eqnarray}\label{e1.2}
 N\{\gamma\} = \sup_{u \in \F}\ \inf_{T \in \gamma} |T(du \wedge d^c u \wedge (dd^cu)^{k-1})|, \quad \mbox{ if } \dim\gamma = 2k,
\end{eqnarray}
where $T$ runs over all currents which represent $\gamma$ and 
$${\mathcal F} = \{u \in C^2(M)\mid u \mbox{ is plurisubharmonic and $0 < u < 1$ on $M$}\}.$$

The intrinsic norm $N$ may also be obtained as the supremum over the subclass 
of $C^2$ solutions of homogeneous complex
Hessian equations,
\begin{eqnarray*}
\F'_k = \{u \in \F\mid (dd^c u)^k = 0, \dim \gamma = 2k-1, \mbox{ or } du \wedge (dd^c u)^k = 0, \dim\gamma = 2k\}.
\end{eqnarray*}
When $k = 2n-1$, elements of $\F'_{2n-1}$ are plurisubharmonic functions satisfying the homogeneous complex Monge-Amp\`ere equation
\begin{eqnarray}\label{e1.4}
(dd^cu)^n = 0. 
\end{eqnarray}

In a series of works \cite{BT1}--[7], based on the fundamental Chern-Levine-Nirenberg inequality and Lelong's work on positive currents, 
Bedford-Taylor developed weak solution theory for complex Monge-Amp\`ere equations. They extended Monge-Amp\`ere operators for locally bounded plurisubharmonic functions as positive currents, and obtained many
important results. In \cite{BT3} and \cite{BT4}, the variational properties
of complex Monge-Amp\`ere equations were studied. In particular in \cite{BT4},
Bedford-Taylor introduced the extended intrinsic norm $\tilde N$.
\begin{eqnarray}\label{e1.9} \qquad
\tilde N\{\gamma\} &\hskip-7pt =\hskip-7pt& \sup_{u \in \tilde{\F}}\ \inf_{T \in \gamma} |T(d^cu \wedge (dd^cu)^{k-1})|, 
\quad \mbox{ if } \dim \gamma = 2k-1,\\
\tilde N\{\gamma\} &\hskip-7pt =\hskip-7pt& \sup_{u \in \tilde{\F}}\ \inf_{T \in \gamma} |T(du \wedge d^cu \wedge (dd^cu)^{k-1})|, 
\quad \mbox{ if } \dim \gamma = 2k, \label{e1.10}
\end{eqnarray}
where the infimum this time is taken over smooth, compactly supported currents which represent $\gamma$ and
$$\tilde{\F} = \{u \in C(M)| u \mbox{ is plurisubharmonic }, 0 < u < 1 \mbox{ on } M\}.$$

It is shown that $\tilde N \geq N$ and $\tilde N
<\infty$. $\tilde N$ also enjoys other similar properties of $N$. They are
invariants of the complex structure, and decrease under holomorphic maps.  These properties are useful in the study of holomorphic mappings.
The intrinsic norms can be extended to intrinsic pseudo-metrics on the manifold which are closely related to
Caratheodory and Kobayashi metrics. 

It was observed in \cite{CLN} that equation (\ref{e1.4}) also arises as the Euler equation for a stationary point of the convex functional
\begin{eqnarray}\label{e1.6}
I(u) = \int_M du \wedge d^cu \wedge (dd^cu)^{n-1}.
\end{eqnarray}

Let $M$ be a closed complex manifold with smooth boundary $\partial M=\Gamma_1 \cup \Gamma_0$, and let
\begin{eqnarray}\label{e1.200}
\B = \{ u \in \F \mid u = 1 \mbox{ on $\Gamma_1$,  $u=0$ on $\Gamma_0$}\}.
\end{eqnarray}
If $v \in \B$, let $\gamma$ denote the $(2n-1)$-dimensional homology class of the level hypersurface $v = $ constant.  Then
for all $T \in \gamma$, if $v$ satisfies $(dd^cv)^n = 0$,
$$\int_T dv \wedge (dd^cv)^{n-1} = \int_M dv \wedge d^cv \wedge (dd^cv)^{n-1} = I(v).$$ 
\pagebreak

 The following conjecture was made in \cite{CLN} by Chern-Levine-Nirenberg.
 
\nonumproclaim{{C}onjecture} $N\{\Gamma_1\} = \inf_{u \in \B} I(u).$
\endproclaim
 
Suppose $M$ is of the following form,
\begin{eqnarray}\label{e1.8}
M = \bar \Omega^* \setminus \left(\bigcup ^N_{j=1} \Omega_j\right),
\end{eqnarray}
where $\Omega^*$, $\Omega_1,\ldots,\Omega_N$ are bounded strongly smooth pseudoconvex domains in ${\bf C}^n$
where 
$\bar\Omega_1,\ldots,\bar\Omega_N$ are pairwise disjoint, and $\bar\Omega_j \subset \Omega^*$, for all $j =
1,\ldots,N$.  Also,  $\bigcup ^N_{j=1} \Omega_j$ is holomorphically convex in
$\Omega^*$ (this is a necessary condition in order for equation (\ref{e1.12})
to have a solution), and $\Gamma_1=\partial \Omega^*$ and $\Gamma_0=\bigcup ^N_{j=1} \partial \Omega_j$.
For the solution $u$ of (\ref{e1.12}), Bedford-Taylor proved in \cite{BT4} that
\begin{eqnarray}\label{e1.101}
\tilde N(\{\Gamma_1\}) = \int_M du \wedge d^cu \wedge (dd^cu)^{n-1},
\end{eqnarray}
and 
\begin{eqnarray}\label{e1.102}
\tilde N(\{\Gamma_1\}) = \inf_{v \in \tilde {\cal B}} \int_M dv \wedge d^cv \wedge (dd^cv)^{n-1},
\end{eqnarray}
where
$$\tilde {\cal B} = \{u \in \tilde{\F} | u \in {\rm Lip}(M) \cup C^2(\partial M), \ u \geq 1 \mbox{ on } \Gamma_1,\ u \leq 0 \mbox{ on }
\Gamma_0\}.$$ Furthermore, if the extremal function is in $C^2(M)$, the identity (\ref{e1.20}) is valid. And if $u \in C^1(M)$, then 
$$\tilde N(\{\Gamma_1\})\leq \int_{\Gamma_1} \left(\frac{\partial u}{dr}\right)^n d^cr \wedge (dd^cr)^{n-1}.$$

Concerning the manifold $M$ in the conjecture, we
remark that if $\Gamma = \{ v = \mbox{constant}\}$ for some $v \in \B$, $\Gamma \sim \{v=1\} \sim \{v=0\}$ in $H_{2n-1}(M)$.  The
hypersurface $\{v=1
\}$ is pseudoconvex, and $\{v=0\}$ is pseudoconcave.  If $M$ is embedded in
${\bf C}^n$, $v$ is strictly plurisubharmonic, and $M$ must be of the form
(\ref{e1.8}).

In fact, we will show the reverse is also true. The following proposition will be used in the estimation of the degenerate Monge-Amp\`ere equations via the subsolution method.

\proclaim{Proposition}\label{subsol}
If $M$ is of the form {\rm (\ref{e1.8}),} there is $v \in {\rm PSH}(M^0) \cap C^\infty(M)$ such that 
$v = 1+cr$ near $\Gamma_1$ and $v = c_jr_j$ near $\partial\Omega_j${\rm ,} $j = 1,2,\ldots,N${\rm ,}
 for some positive
constants $c${\rm ,} $c_1,\ldots,c_N${\rm ,}
 where $r$ and $r_j$ are the defining functions of $\Omega$ and $\Omega_j${\rm ,} $j=1,
\ldots, N${\rm ,} respectively{\rm .}  Furthermore{\rm ,}
\begin{eqnarray}\label{e12.3}
(dd^cv)^n > 0 \qquad \mbox{ in } M.
\end{eqnarray}
\endproclaim

Therefore, we will concentrate on the manifolds of the form  (\ref{e1.8}).
We now state our main results.

\proclaimtitle{Regularity Theorem}
\specialnumber{1.1} \proclaim{Theorem}  If $M$ is of the form {\rm (\ref{e1.8}),} for the
unique solution
$u$ of {\rm (\ref{e1.12}),} there is a sequence $\{u_k\} \subset \B$ such that
$$\|u_k\|_{C^2(M)} \leq C, \quad \hbox{for all } \ k, \qquad \lim_{k \to \infty} \sup(dd^cu_k)^n = 0,$$
and $\lim_{k \to \infty} \|u_k - u \|_{C^{1,\alpha}(M)} = 0, \quad \hbox{for all }\ 0 < \alpha < 1${\rm
.} In particular{\rm ,}
 $u \in C^{1,1}(M).$ 
\endproclaim

As a consequence of the regularity theorem, we establish identity (\ref{e1.20}).
This identity can be used to obtain the uniqueness of the extremal functions (e.g., Lemma 2.6 in \cite{BB}) and, in turn, to study holomorphic mappings. We will discuss these applications elsewhere.

\proclaim{Theorem} If $M$ is of the form {\rm (\ref{e1.8}),} 
\begin{eqnarray}\label{e1.21}
\tilde N(\{\Gamma_1\})
=\int_{\Gamma_1} \left(\frac{\partial u}{dr}\right)^n d^cr \wedge (dd^cr)^{n-1}
\end{eqnarray}
where $r$ is the defining function of $\Omega.$
\endproclaim

The Chern-Levine-Nirenberg conjecture for $\tilde N$ is also valid. That is,
the class $\tilde {\cal B}$ in (\ref{e1.102}) can be replaced by a smaller class
${\cal B}$ defined in (\ref{e1.200}).

\proclaim{Theorem} If $M$ is of the form {\rm (\ref{e1.8}),} then 
\begin{eqnarray}\label{e1.22}
\tilde N(\{\Gamma_1\}) = \inf_{v \in {\cal B}} \int_M dv \wedge d^cv \wedge (dd^cv)^{n-1}.
\end{eqnarray}
\endproclaim

The rest of the paper is organized as follows. In Section 2, we establish
$C^2$ {\it a priori} estimates for the solutions of the equation under the
assumption of the existence of a strictly plurisubharmonic solution. 
Section 3 is devoted to the construction of the subsolution and the proofs
of the main results.

\demo{Acknowledgement}  This work was completed while the
author was visiting NCTS at National Tsing Hua University in
Taiwan. He would like to
thank Professor C.S. Lin for the kind arrangement and thank
NCTS for its warm hospitality.\enddemo

 \section{Estimates for the  extremal function}

Equation (\ref{e1.12}) is a degenerate complex Monge-Amp\`ere equation.
There have been many works in this direction. In particular, the work of
Caffarelli-Kohn-Nirenberg-Spruck \cite{CKNS} establishes $C^{1,1}$
regularity
for solutions in strongly pseudoconvex domains with homogeneous boundary
condition. On the other hand, the $C^{1,1}$ regularity of the homogeneous
complex Monge-Amp\`ere equation on strongly pseudoconvex domains with
arbitrary boundary data was obtained by Krylov \cite{Kr}.
There is also extensive literature on degenerate real Monge-Amp\`ere
equations (see \cite{PG} and \cite{GTW} for the references).
In our case, some pieces of the boundary are concave. This is the
main technical difficulty of the problem. We will incorporate 
the subsolution method of \cite{G} where interior regularity was 
treated for pluricomplex Green functions.  
The existence of a smooth subsolution $v$ of (\ref{e2.1}) in 
Proposition~\ref{subsol} is crucial in the proof of the theorem. 
This proposition will be proved in the next section.

To establish the existence of the solution $u \in C^{1,1}(M)$, we consider
the following equation with parameter $ 0\leq t <1$,
\begin{eqnarray}\label{e2.1}
\left\{\begin{array}{l}
(dd^cu)^n = (1-t)f_0 \\
u|_{\Gamma_1} = 1 \\
u|_{\Gamma_0} = 0 ,
\end{array}\right.
\end{eqnarray}
where $f_0= \det(v_{i \bar j})$ and $v$ is as in Proposition~\ref{subsol}. When
$0 \leq t <1$, the equation is elliptic. We want to prove that equation 
(\ref{e2.1}) has a unique smooth solution
with a uniform $C^{1,1}$ bound. 

\proclaim{Theorem}
With $M$ as in {\rm (\ref{e1.8}),} there is a constant $C$ depending only
on $M$ {\rm (}\/independent of $t${\rm )}
such that for each $0\leq t <1${\rm ,} there is a unique smooth solution $u$ of {\rm (\ref{e2.1})} with 
\begin{eqnarray}\label{e2.2}
\|u\|_{C^{1,1}(M)} \leq C .
\end{eqnarray}
\endproclaim

The {\it a priori} estimate (\ref{e2.2}) is the key to the
proof of the main results. $C^{1,1}$ regularity of   equation (\ref{e1.12})
can be deduced directly from the above theorem. In the following proof,
we indicate $c$ to be the constant (which may vary line
to line) depending only on $M$ (independent of $t$). We assume the
existence of a subsolution $v$ in Proposition~\ref{subsol}.

\demo{{P}roof of Theorem {\rm 2.1}} Let $f_0 = \det(v_{i\bar j})$. Then $f_0 > 0$, $f_0 \in C^\infty(M)$.
We consider $0 \leq t \leq 1$.
\begin{eqnarray}\label{e2.4}
\left\{\begin{array}{l}
\det(u_{i\bar j}) = (1-t)f_0\\
u|_{\Gamma_1} = 1\\
u|_{\Gamma_0} = 0.
\end{array}\right.
\end{eqnarray}
We show that for all $0 \leq t < 1$, there exists $u_t \in C^\infty$, $u_t$ plurisubharmonic, such that $u_t$ solves
(\ref{e2.4}) and there exists $C > 0$, for all $0 \leq t < 1$
\begin{eqnarray}\label{e2.5}
\|u_t\|_{C^{1,1}(M)} \leq C.
\end{eqnarray}
The uniqueness is a consequence of the comparison theorem for complex Monge-Amp\`ere equations (see \cite{BT1},  \cite{BT2} and \cite{CKNS}).  We also note that if (\ref{e2.5}) holds, the equation is elliptic when $0 \leq t < 1$.  By
the  Krylov-Evans theorem, $u_t \in C^\infty(\bar M)$.  Therefore, everything is reduced to the establishment of the {\it a
priori} estimates (\ref{e2.5}).  In the rest of the proof, we will drop the subindex $t$.  
 \enddemo

 $C^0$-{\it estimates}. Since $u$ is plurisubharmonic in $M^0$, and $0 \leq u \leq 1$ on $\partial M$, the maximum
principle gives $0 \leq u(z) \leq 1$ for all $z \in M$.

\vglue6pt
 $C^1$-{\it estimates}. For any first order differential operator $D$ with constant real coefficients, we consider $w =
Du+Av$, where $A$ is a constant to be picked up later.  Apply $D$ to the equation (\ref{e2.1}); we get
$$u^{i\bar j}(Du)_{i\bar j} = \frac{Df_0}{f_0}\ .$$

We also have
$$u^{i\bar j}(Av)_{ij} \geq Ac \sum u^{i\bar i} \geq n Acf_0^{-\frac{1}{n}}(t-1)^{-\frac{1}{n}} \geq nAcf_0^{-\frac{1}{n}}.$$
If $A$ is large,  
$$u^{i\bar j} w_{i\bar j} > 0 \qquad \mbox{ in } M.$$
By the Maximum Principle for elliptic operators,
$$\max_{\Omega}(Du + Av) = \max_{\partial M}(Du + Av).$$
To estimate $Du$ on $\partial M$, let  $h$ be the solution of
\begin{eqnarray}\label{e2.100}
\left\{\begin{array}{ll}
\Delta h = 0 & \mbox{ in } M^0\\
h|_{\Gamma_1} = 1\\
h|_{\Gamma_0} = 0.
\end{array}\right.
\end{eqnarray}
Since $\det(u_{i\bar j}) \leq \det (v_{i \bar j})$,  
$$\Delta u = \sum^n_{i=1} u_{i \bar i} = nf^{-\frac{1}{n}}_0 (1-t)^{\frac{1}{n}} > 0 = \Delta h $$
and
$$u|_{\partial M} = v|_{\partial M} = h|_{\partial M},$$
by the Comparison Principle, $v(z) \leq u(z) \leq h(z)$, for all $z \in M$.  Therefore
\begin{eqnarray}\label{e2.6}
|Du(z)| \leq \max\left(|Dv(z)|,|Dh(z)|\right) \leq c \qquad \hbox{for all } \ z \in \partial M.
\end{eqnarray}
That is, $\max_{\partial M}|Du| \leq c$.  In turn,
\begin{eqnarray}\label{e2.7}
\max_M W \leq c, \qquad \mbox{ and } \qquad \max_M |Du| \leq c.
\end{eqnarray}

\vglue4pt $C^2$-{\it estimates}. Since $\det^{\frac{1}{n}}$ is concave,  
$$u^{i\bar j}(D^2u)_{i \bar j} \geq n \ \frac{D^2((f)^{\frac{1}{n}})}{(f)^{\frac{1}{n}}} = n\ \frac{D^2f^{\frac{1}{n}}_0}{f^{\frac{1}{n}}_0}\ .$$
This yields 
\begin{eqnarray*}
u^{i\bar j}((D^2u + Av)_{i \bar j}) &\geq& n   \frac{D^2f^{\frac{1}{n}}_0}{f^{\frac{1}{n}}_0} + Ac \Sigma u^{i \bar i}\\ 
&\geq& n\ \frac{D^2f^{\frac{1}{n}}_0}{f^{\frac{1}{n}}_0} + nAcf^{-\frac{1}{n}}_0\\
&\geq& n\ \frac{D^2f^{\frac{1}{n}}_0}{f^{\frac{1}{n}}_0} + nAcf^{-\frac{1}{n}}_0.
\end{eqnarray*}
If we pick $A \geq \frac{1}{c} \max_{z \in M} |D^2f^{\frac{1}{n}}_0|$, then
$$u^{i\bar j} ((D^2u + Av)_{i\bar j}) \geq 0 \qquad \mbox{ in } \Omega.$$
By the Maximum Principle,
\begin{eqnarray}\label{e2.8}
\max_M\{D^2u + Av\} = \max_{\partial M}\{D^2u + Av\}.
\end{eqnarray}

In order to obtain $C^2$-{\it a priori} estimates for $u$, we need to get the estimates of the second derivatives of $u$ on the boundary of
$M$.  The boundary of $M$ consists of pieces of compact strongly pseudoconvex and pseudoconcave hypersurfaces. 
The second derivative estimates on strongly pseudoconvex hypersurface
were established in \cite{CKNS}. Therefore, we will concentrate
on the pseudoconcave piece $\Gamma_0$ of the boundary. 
We follow the arguments in 
\cite{G} with the aid of a strictly plurisubharmonic subsolution $v$ of Proposition~\ref{subsol}. Let $h$ be the 
harmonic function of (\ref{e2.100}) in $M$. 

Suppose $z_0 \in \Gamma_0$, let $z_1 = x_1 + \sqrt{-1}y_1,\ldots, z_{n-1} = x_{n-1} + \sqrt{-1} y_{n-1}$, and $z
_n = x_n + \sqrt{-1}y_n$.  We may assume $z_0 = 0$ and that
$\frac{\partial}{\partial z_j}$,   $\frac{\partial}{\partial \bar z_j}$, $\frac{\partial}{\partial y_n}$ are tangential to $\partial M$ at
$z_0$.  Set $t_1 = x_1$, $t_2  = y_1, \ldots,t_{2n-3} = x_{n-1}$, $t_{2n-2} = y-{n-1}$, $t_{2n-1} = y_n = t$.  
We also assume $\frac{\partial v}{\partial x_n}=-1$ at $z_0=0$.
\vglue2pt
For all $\varepsilon >0$, we let $S_\varepsilon = B_\varepsilon(0) \cap M^0$, 
and   define,
\begin{eqnarray}\label{e2.10}
w(z) = (v(z)-u(z)) +a(v(z)-h(z)) + bv^2(z).
\end{eqnarray}
 The constants $a$ and $b$ will be chosen later.
The following lemma was proved by B. Guan (Lemma 2.1 in \cite{G}). For  completeness, we reproduce the proof (with some minor modification) here.

\specialnumber{2.1} \proclaim{Lemma}
Let $\tilde c=\inf_{z\in M, \xi\in {\bf C^n}\setminus \{0\}}\frac{\sum v_{i\bar j}\xi_i\bar \xi_j}{|\xi|^2}${\rm .}
For suitable choices of positive constants $a, b$ and $\varepsilon${\rm ,}
\begin{eqnarray}
\mbox{\rm (i)} \qquad Lw(z) &\geq& \frac{\tilde c}{4}\sum u^{i\bar i}(z)  \qquad \mbox{ in } S_\varepsilon,\label{e2.11}\\
&&\nonumber\\
\mbox{\rm (ii)} \qquad  w(z)& \leq &0\phantom{\sum u^{i\bar i}(z) \, } \qquad \mbox{ on }  S_\varepsilon, \label{e2.12}
\end{eqnarray}
where 
\begin{eqnarray}\label{e2.200}
L = \sum^n_{i,j} u^{i \bar j} \frac{\partial^2}{\partial z_i \partial \bar z_j}.
\end{eqnarray} 
\endproclaim

{\it Proof of Lemma} 2.1.
We first observe that,
\begin{eqnarray}\label{e2.201}
L(v-u)\geq Lv-n.
\end{eqnarray}
Also, as ${\rm tr}(v_{i \bar j}) \geq n \tilde c$ (and by Hopf's lemma),
\begin{eqnarray}\label{e2.212}
h(z)-v(z) \geq c_0 v(z), \quad \hbox{for all }\  z \in S_\varepsilon,
\end{eqnarray}
with some uniform constant $c_0>0$. Furthermore,
$$L(v-h)\geq -C_1\sum u^{i\bar i},$$
for some uniform positive constant $C_1$ as $h \in C^2$. 
We also have,
$$L(v^2)=2vLv+2 \sum u^{i\bar j}v_iv_{\bar j}\geq 0.$$
Therefore, as $a$ and $b$ are positive,
\begin{eqnarray}\label{e2.205}
Lw \geq -aC_1\sum u^{i\bar i}+(1+2bV)Lv +2b \sum u^{i\bar j}v_iv_{\bar j}-n.
\end{eqnarray}
Now, as $\frac{\partial v}{\partial x_n}=-1$ at $z_0=0$ and $Tv(0)=0$ for any
tangential vector field $T$, if $\varepsilon>0$ is small,
\begin{eqnarray}\label{e2.203}
\sum u^{i\bar j}v_iv_{\bar j}\geq u^{n\bar n}v_nv_{\bar n}-C_3\varepsilon \sum u^{i\bar i}
\end{eqnarray}
for some positive constant $C_3$ under control. 
By the geometric inequality,
\begin{eqnarray}\label{e2.206}
\frac{\tilde c}{4}\sum u^{i\bar i} + bu^{n\bar n}\geq \frac{n(\frac{\tilde c}{4})^{\frac{n-1}{n}}(\frac{\tilde c}{4}+b)^{\frac1n}}
{(\prod_i u_{i\bar i})^\frac1n} \geq \frac{n(\frac{\tilde c}{4})^{\frac{n-1}{n}}b^{\frac1n}}
{\det^{\frac1n}u_{i\bar j}}.
\end{eqnarray}
Since $\det(u_{i\bar j})=(1-t)f_0 \leq f_0$, we pick $b$ such that 
\begin{eqnarray}\label{e2.210}
\inf_{z \in M} \{\frac{n(\frac{\tilde c}{4})^{\frac{n-1}{n}}b^{\frac1n}}
{\det^{\frac1n}v_{i\bar j}(z)} \}\geq n.
\end{eqnarray}
Putting (\ref{e2.206}), (\ref{e2.203}) and (\ref{e2.210}) into (\ref{e2.205}), we obtain 
\begin{eqnarray}\label{e2.211}
Lw \geq -(aC_1+2C_3b\varepsilon+\frac{\tilde c}{4}\sum u^{i\bar i})+(1+2bv)Lv. 
\end{eqnarray}
Since $v>0$ is dominated by $\varepsilon$ in $S_{\varepsilon}$, and 
$Lv\geq \tilde c \sum u^{i\bar i}$, if we choose $a=\frac{\tilde c}{1+C_1}$ and $\varepsilon$ small, (\ref{e2.11}) follows from (\ref{e2.211}).
\pagebreak
 
To examine values of $w$ in $S_\varepsilon$, we note that
$v=0$ on $\partial \Gamma_0$. By Hopf's lemma, if $\varepsilon$ is sufficiently small, for $c_0$ as in (\ref{e2.212}),
$$w\leq -ac_0v+bv^2\leq 0.$$ 
The proof of the lemma is complete. \hfill \qed

\bigskip
Now, back to the proof of Theorem 2.1: 
Let
\begin{eqnarray}\label{e2.9}
T_\alpha = \frac{\partial}{\partial t_\alpha} - \frac{v_{t_\alpha}}{v_{x_n}}\ \frac{\partial}{\partial x_n}.
\end{eqnarray}
Let
\begin{eqnarray*}
\tilde w(z)= w(z) = \pm T_\alpha u(z) + u^2_t(z) +\tilde C w(z)-B|z|^2,
\end{eqnarray*}
where $B$ and $\tilde C$ are certain constants to be picked later.

By \cite{CKNS},
\begin{eqnarray*}
L(\pm T_\alpha u(z) + u^2_t(z)) \geq -c(1+\sum u^{\i\bar i}),
\end{eqnarray*}
in $S_\varepsilon$. From the above lemma, if we
pick $\tilde C$ and $B$ large enough, with $\tilde C$ sufficiently larger than $B$ (as
$\sum u^{\i\bar i}\geq \det^{-\frac1n}u_{i\bar j}\geq f_0$), we will have
\begin{eqnarray}
\mbox{\phantom{i}(i)} && L(\tilde w)(z) \geq 0 \qquad \mbox{ in} S_\varepsilon,\label{e2.111}\\
&&\nonumber\\
\mbox{(ii)} && \tilde w(z) \leq 0\phantom{\, \tilde w)} \qquad \mbox{ on } \partial S_\varepsilon, \label{e2.122}
\end{eqnarray}

Since $w(0) = 0$, we must have $w_{x_n}(0) \leq 0$.  As a consequence, $|u_{x_nt_\alpha}(0)| \break \leq C$. 

Near $\Gamma_0$,
 $u(z) = \sigma(z)v(z)$. But $u(z) \geq v(z) \geq 0$, for all $z \in M$, so  that $\sigma(z) \geq 1$. On the other hand, if $\vec n$ is the
normal to
$\Gamma_0$,
$$u_n=\sigma v_n,  \quad \hbox{for all }\  z \in \Gamma_0.$$ 
As $|Du|$ is bounded in $M$, $u_n\leq 0$ and $v_n=-1$ on $\Gamma_0$, 
so $1 \leq \sigma \leq c$.
On $\Gamma_0$, we have
$$u_{i\bar j}(z) = \sigma(z)v_{i \bar j}(z), \qquad \hbox{for all } \ i,j \leq n-1.
$$
It follows that 
\begin{equation}\label{tang}
cI_{n-1} \geq u_{i\bar j}(z) \geq v_{i \bar j}(z) \geq I_{n-1}, \quad \hbox{for all } \ z \in \Gamma_0, \quad \hbox{for all } \ i,j \leq n-1,
\end{equation}
where $I_{n-1}$ is the $(n-1)$-dim identity matrix.

This gives tangential second derivative bounds.
To estimate $u_{x_nx_n}(0)$, we note that  
$$u_{x_nx_n}(0) = 4u_{n\bar n}(0) - u_{y_ny_n}(0).$$
We only need to estimate $u_{n \bar n}(0)$.  Since $u_{i\bar n}(0)$ and $u_{n \bar i}(0)$ are bounded, $u_{n \bar n}(0)$ can be estimated by
(\ref{tang})   and equation (\ref{e2.4}).  In conclusion,
\begin{eqnarray}\label{e2.17}
\|u\|_{C^2(M)} \leq c, \qquad \hbox{for all } \ 0 \leq t < 1.
\end{eqnarray}

Finally, the higher regularity for the solution of the equation (\ref{e2.4}) for $0\leq t<1$ follows from the Evans-Krylov theorem and the
Schauder theorem. The existence follows from standard elliptic theory for the Dirichlet problem. 
\phantom{mwsw}
\hfill\qed

\section{Existence of subsolutions}

We construct
 in this section a strictly smooth plurisubharmonic function with the same boundary value as the solution $u$ of equation (\ref{e1.12}). We begin with some useful elementary lemmas.

\specialnumber{3.1} \proclaim{Lemma}\label{lem1}
For all $\delta > 0${\rm ,} there is an even function $h(t) \in C^\infty({\bf R})${\rm ,} such that
\begin{itemize}
\ritem{(i)} $h(t) \geq |t|${\rm ,}  for all $t \in {\bf R}${\rm ,}x $h(t) = |t|${\rm ,}  for all $|t| \geq \delta${\rm ;}

\ritem{(ii)} $|h'(t)| \leq 1$ and $h^{\prime\prime}(t) \geq 0${\rm ,}   for all $t \in {\bf R}$ and $h'(t) \geq 0${\rm ,} for all $t \geq 0${\rm
.}
\end{itemize}

\endproclaim

{\it Proof}. For all $\varepsilon > 0$, let $\rho_\varepsilon(s) = \frac{1}{\varepsilon}\rho\left(\frac{s}{\varepsilon}\right)$, where
$\rho$ is a standard nonnegative smooth even function supported in $|s| \leq 1$ with $\int^\infty_{-\infty} \rho(s)ds = 1$.  Set
$$\tilde h_\varepsilon (t) = \int^t_{-\infty} \rho_\varepsilon(s)ds - \int^\infty_t \rho_\varepsilon(s)ds.$$
Now,
\begin{eqnarray*}
\tilde h'_\varepsilon(t) &=& 2\rho_\varepsilon(t) \geq 0,\\
|\tilde h_\varepsilon(t)| & \leq & \int^t_{-\infty} \rho_\varepsilon(s)ds + \int^\infty_t \rho_\varepsilon(s)ds\\
&=& \int^\infty_{-\infty} \rho_\varepsilon(s)ds = 1, \qquad \hbox{for all } \ t \in {\bf R}.
\end{eqnarray*}
Moreover,
$$\tilde h_\varepsilon(t) = 1,  \mbox{ if } t > 2\varepsilon, \quad \hbox{and}
\quad 
\tilde h_\varepsilon(t) = -1, \mbox{ if } t < -2\varepsilon.$$
We also have,
$$\tilde h_\varepsilon(t) \geq 0,  \mbox{ if } t \geq 0,\quad \hbox{and}
\quad \tilde h_\varepsilon(t) \leq 0,  \mbox{ if } t \leq 0.$$
Since $\rho_\varepsilon(s)$ is nonnegative and even, 
$$\tilde h_\varepsilon(t) = -\tilde h_\varepsilon(-t).$$
Now, for $\varepsilon = \frac{\delta}{2}$, set
$$h(t) = \int^t_0 \tilde h_\varepsilon(s)ds + \left( 1 - \int^1_0 \tilde h_\varepsilon(s)ds\right).$$
We have
$$h(t) = h(-t), \qquad h'(t) = \tilde h_\varepsilon(t),$$
so that $th'(t) \geq 0$,   for all $t$.  Also 
$$|h'(t)| = |\tilde h_\varepsilon(t)| \leq 1.$$
Since $h(1) = 1$ and $h'(t) =1$,   for all $t \geq \delta$,  we have $ h(t) = t$,  for all $t \geq \delta$. Therefore,  $ h(t)= |t|$,   for all $|t|
\geq \delta$, as $h$ is even. Finally
\vglue4pt 
\hfill $h^{\prime\prime}(t) = \tilde h'_\varepsilon(t) = 2\rho_\varepsilon(t) \geq 0 \qquad \hbox{for all }\ t \in {\bf R}.$ 
\hfill\qed

\specialnumber{3.2} \proclaim{Lemma}\label{lem2}
Suppose $\Omega$ is a domain of $\bf C^n${\rm .}
 For $f,g \in C^k(\Omega)${\rm ,} $k \geq 2${\rm ,}  for all $\delta > 0${\rm ,} there is an $H \in C^k(\Omega)$ such that
\begin{itemize}
\ritem{(i)} $H \geq \max(f,g)$ and
$$ H(z) = \left\{\begin{array}{ll}
f(z), & \mbox{ if } f(z) - g(z) > \delta\\
g(z), & \mbox{ if } g(z) - f(z) > \delta;
\end{array}\right.$$

\ritem{(ii)} There exists $|t(z)| \leq 1${\rm ,} such that
$$(H_{ij}(z)) \geq \left(\frac{1+t(z)}{2} f_{i\bar j}(z) + \frac{1-t(z)}{2} g_{ij}(z)\right), \ \hbox{for all } \ z \in \{|f-g| < \delta\}.$$
\end{itemize}
\endproclaim

{\it Proof}. Note that $\max(f,g) = \frac{f+g}{2} + \frac{|f-g|}{2}$.  Let $h$ be the function as in Lemma \ref{lem1}.  We set
$$H(z) = \frac{f(z) + g(z)}{2} + \frac{h(f(z) - g(z))}{2}
.$$
It is obvious that $H$ satisfies property (i).  Now we calculate $H_{ij}(z)$: \begin{eqnarray*}
H_{i\bar j}(z) &=& \frac{f_{i\bar j}(z) + g_{i\bar j}(z)}{2} + \frac{1}{2} \biggl\{h'(f(z) - g(z))(f_{i\bar j}(z) - g_{i \bar j}(z))\\
&& +\ h^{\prime\prime}(f(z) - g(z))(f(z) - g(z))_i(f(z) - g(z))_{\bar j}\biggr\}\\
&\geq& \frac{1+ h'(f(z) - g(z))}{2} f_{i\bar j}(z) + \frac{1- h'(f(z) - g(z))}{2} g_{i\bar j}(z).\\
\noalign{\vskip-36pt}
\end{eqnarray*} 
\phantom{ouch}\hfill\qed
\vglue24pt

We now construct a subsolution $v$ in   Proposition 2.1.
 
\demo{Proof of Proposition {\rm 2.1}}  A Lipschitz continuous subsolution can be
constructed as in \cite{BT4}.  Let $\psi_0,\ldots,\psi_N$ be the 
defining functions of $\Omega^*, \Omega_1,\ldots,\Omega_N$ respectively, such that $\psi_j\in C^\infty(\bar\Omega_0)$ and $dd^c\psi_j > 0$ in a neighborhood $U_j$ of $\partial\Omega_j$,
for all $j = 1,\ldots,N$.  (Note that $\psi_j < 0$ in $\Omega_j$ and $\psi_j \neq 0$ on $\partial\Omega_j$.)  Since $\bigcup^N_{j=1}
\Omega_j$ are holomorphic convex in $\Omega_0$, there is a plurisubharmonic function $\psi$ in $\Omega_0$, such that $\psi(z) < 0$ in a
small neighborhood of $\bigcup^N_{j=1} \partial \Omega_j$ and is  positive outside of $\bigcup^N_{j=1} \bar U_j$.  We may assume that
$U_j$ are pair-wise disjoint for $j=1, \ldots, N$, and $\psi \in C^\infty(\bar \Omega_0)$ and $dd^c\psi > 0$ in $\bar\Omega_0$. In each
$U_j$, we set
$$V_j(z) = \max\{\varepsilon ^2\psi_j, \varepsilon\psi
\}.$$    
If we pick $\varepsilon$ small, we have $V_j(z) = \varepsilon^2\psi_j(z)$ in a small neighborhood of $\partial \Omega_j$ in $U_j$, and $V_j(z) = \varepsilon\psi(z)$ away from that small neighborhood in $U_j$. So, we may
extend $V_j$ to all of $\Omega$ by defining 
$V_j(z) = \varepsilon\psi(z)$ outside of $U_j$ and $V_j(z) \leq 1$
on $\Gamma_1$. By applying  Lemma \ref{lem2} in $U_j$, we obtain a smooth strictly plurisubharmonic function $H_j$ such that $H_j=0$ on $\partial \Omega_j$ and $H_j = \varepsilon\psi(z)$ outside of $U_j$. Now, we pick $\lambda$ large so that $1 + \lambda \psi_0$ will be
very negative in $\bigcup^N_{j=1} \bar U_j$.  Set 
$$H(z) = \max\{1 + \lambda \psi_0(z), H_1(z),\ldots, H_N(z)\}.$$  

The function $H$ is taken as a maximum of smooth strongly  
plurisubharmonic functions, and at no point are there  three equal 
functions involved. Therefore, we can apply Lemma \ref{lem2} near 
the places where 
any two of these function meet to obtain a smooth strongly plurisubharmonic function $v$, such that $v(z) = 1$ on $\Gamma_1$, $v(z) = 0$ on $\Gamma_0$ and $(dd^c)^nv > 0$ in $M^0$.
\enddemo

\demo{{P}roof of the Regularity Theorem {\rm 1.1}} From the above proposition and Theorem 2.1, there is a sequence of strictly smooth
plurisubharmonic functions $\{u^t\}$ satisfying (\ref{e2.1}). By (\ref{e2.2}), there is a subsequence $\{t_k\}$ that tends
to $1$, such that $\{u_{t_k}\}$ converges to a plurisubharmonic function $u$
in $C^{1, \alpha}(M)$ for any $0<\alpha<1$. By the
Convergence Theorem for complex
Monge-Amp\`ere measures, $u$ satisfies equation (\ref{e1.12}). Again by
(\ref{e2.2}), $u \in C^{1,1}(M)$. \hfill\qed\enddemo

\demo{Proof of Theorem {\rm 1.2}} If we let $\{u_k\}$ be as in the Regularity Theorem, we have
\begin{eqnarray*}
\int_M du_k \wedge d^cu_k \wedge (dd^cu_k)^{n-1} &=&
\int_{\Gamma_1} d^cu_k \wedge (dd^cu_k)^{n-1} - \int_M u_k(dd^cu_k)^n\\
&=& \int_{\Gamma_1} \left(\frac{\partial u_k}{\partial r}\right)^n d^c r \wedge (dd^c r)^{n-1} - \int_M u_k(dd^cu_k)^n.
\end{eqnarray*}
Since $u_k \to u$ in $C^{1,\alpha}(M)$, $\left(\frac{\partial u_k}{\partial r}\right)^n \to 
\left(\frac{\partial u}{\partial r}\right)^n$ uniformly on $\Gamma_1$.  Therefore, by the Convergent Theorem for complex Monge-Amp\`ere measures (see \cite{BT2}), we get
\vglue4pt
\hfill ${\displaystyle \tilde{N}(\{P_1\}) = \int_M  du\wedge d^cu\wedge (dd^cu)^{n-1}= \int_{\Gamma_1} 
\left(\frac{\partial u}{\partial r}\right)^n d^cr \wedge (dd^cr)^{n-1}.}$\enddemo

\demo{Proof of Theorem {\rm 1.3}} By \cite{BT4},  
$$\tilde N(\{\Gamma_1\}) = \int_{\Gamma_1} d^cu \wedge (dd^cu)^{n-1} = \int_M d u \wedge d^cu \wedge (dd^cu)^{n-1},$$
and $\tilde N(\{\Gamma_1\}) \geq N(\{\Gamma_1\})$ by definition.  Also by the
Comparison Theorem (\cite{BT1}, \cite{CKNS}), for all $v \in B$ if $v \not\equiv u $, one must have $ v < u$ in $M^{0}$. By
Theorem 2.1 in~\cite{BT4},
$$\int_M du \wedge d^c u \wedge (dd^cu)^{n-1} \leq
\int_M dv \wedge d^cv \wedge (dd^cv)^{n-1}.$$
That is,
$$\int_M du \wedge d^c u \wedge (dd^cu)^{n-1} \leq \inf_{v \in B} \int_M dv \wedge d^c v \wedge (dd^cv)^{n-1}.$$
On the other hand, by the Convergent Theorem for 
complex Monge-Amp\`ere measures 
$$\liminf_{k \to \infty} \int_M du_k \wedge d^cu_k \wedge (dd^cu_k)^{n-1} = \int_M du \wedge d^c u \wedge (dd^cu)^{n-1}.$$
That is, 
\vglue4pt
\hfill ${\displaystyle\tilde N(\Gamma_1) = \int_M du \wedge d^c u \wedge (dd^cu)^{n-1} = \inf_{v \in B} \int dv \wedge d^c v \wedge
(dd^cv)^{n-1}.}$ 
\enddemo

{\it Remark}. The main results can be generalized to Stein manifolds without
major changes. Estimates $C^1$ and $C^2$  can be obtained from the fact
that there are finite global holomorphic and anti-holomorphic vector fields that generate $T(M)$.
Interior estimates  $C^2$ also follow from Yau \cite{Y}. We can obtain $C^2$ boundary estimates using the same arguments as in this paper,
as they are local estimates.

 \end{document}